\documentclass{article}
\usepackage{amssymb,amsxtra,amsmath,amsthm}

  \newtheorem{proposition}{Proposition}

\begin{document}

\title{On Mordell-Weil groups
of elliptic curves induced by Diophantine triples}
\author{Andrej Dujella}

\markboth{A. Dujella}{Elliptic curves induced by Diophantine
triples}
\pagestyle{myheadings}



\date{}
\maketitle

\begin{center}
{\it

Dedicated to Professor Sibe Marde\v{s}i\'{c} on the occasion of his 80th
birthday.

\vspace{\baselineskip} }
\end{center}

\begin{abstract}
We study the possible structure of the groups of rational points
on elliptic curves of the form $y^2=(ax+1)(bx+1)(cx+1)$, where
$a,b,c$ are non-zero rationals such that the product of any two of
them is one less than a square.
\end{abstract}

 \vspace{1ex}

\section{Introduction}
\label{intro}

Let $E$ be an elliptic curve over $ \mathbb{Q}$.
By Mordell-Weil theorem, the group $E( \mathbb{Q})$ of rational
points on $E$ is a finitely generated abelian group. Hence, it is
the product of the torsion group and $r\geq 0$ copies of infinite
cyclic group:
$$ E( \mathbb{Q}) \simeq E( \mathbb{Q})_{\rm
tors} \times { \mathbb{Z} }^r. $$
By Mazur's theorem, we know that $E( \mathbb{Q})_{\rm tors}$ is
one of the following 15 groups:
$ \mathbb{Z}/n\mathbb{Z}$ with $1\leq n \leq 10$ or $n=12$,
$ \mathbb{Z}/2\mathbb{Z} \times \mathbb{Z}/2m\mathbb{Z}$ with
$1\leq m\leq 4$.
On the other hand, it is not known which values of rank $r$ are
possible. The folklore conjecture is that a rank can be arbitrary large, but it seems to
be very hard to find examples with large rank. The current record
is an example of elliptic curve over $\mathbb{Q}$ with rank $\geq
28$, found by Elkies in May 2006 (see \cite{Elk-28,D-table-records}).

In the present paper, we will study a special case of this problem.
Namely, we will consider only elliptic curves of the form
$ y^2=(ax+1)(bx+1)(cx+1)$, where $\{a,b,c\}$ is a rational Diophantine triple.
Although this is a very special case, it has some relevance for
the more general problem of determining which ranks are possible for
elliptic curves with prescribed torsion group. In particular, we will show
in Section \ref{Z2*Z8} that every elliptic curve over $\mathbb{Q}$ with
torsion group $ \mathbb{Z}/2\mathbb{Z} \times \mathbb{Z}/8\mathbb{Z}$
is induced by some rational Diophantine triple.

A set $\{a_1,a_2,\dots,a_m\}$ of $m$ non-zero integers (rationals)
is called \emph{a (rational) Diophantine $m$-tuple} if $a_i\cdot a_j+1$
is a perfect square for all $1\leq i<j\leq m$.
Diophantus of Alexandria found a rational Diophantine quadruple
$\left\{ \frac{1}{16},\frac{33}{16}, \frac{17}{4}, \frac{105}{16} \right\}$, while the first
Diophantine quadruple in integers, the set $\{1,3,8,120\}$, was found by Fermat
(see \cite{Dick,Dio,Fermat}).
The famous conjecture is that there does not exist a Diophantine quintuple (in non-zero integers)
(see e.g. \cite{Guy,Wal}).
In 1969, Baker and Davenport \cite{B-D} proved that the Fermat's set
$\{1,3,8,120\}$ cannot be extended to a Diophantine quintuple. Recently, it was proved in \cite{D-crelle}
that there does not exist a Diophantine sextuple and there are only finitely many
Diophantine quintuples.

Let $\{a,b,c\}$ be a (rational) Diophantine triple. We define nonnegative rational numbers $r,s,t$ by
$$ ab+1=r^2, \quad ac+1=s^2, \quad bc+1=t^2. $$
In order to extend
this triple to a quadruple, we have to solve the system
\begin{equation} \label{sys}
ax+1=\square,\qquad bx+1=\square,\qquad cx+1=\square.
\end{equation}
It is natural idea to assign to the system (\ref{sys}) the elliptic curve
\begin{equation} \label{ell}
E: \quad y^2=(ax+1)(bx+1)(cx+1).
\end{equation}
Properties of elliptic curves
obtained in this manner and connections between solutions
of the system (\ref{sys}) and the equation (\ref{ell}) were studied in \cite{D-param,D-ja99,D-P2},
but this was mainly for the case when $a,b,c$ are positive integers. In this paper,
we will assume that $a,b,c$ are non-zero rationals, and we will call $\{a,b,c\}$ a
Diophantine triple (hence, omitting the word \emph{rational}).

The coordinate transformation
$x \mapsto \frac{x}{abc}$, $y \mapsto \frac{y}{abc}$
applied on the curve $E$ leads to the elliptic curve
\begin{eqnarray} \label{E'}
 E': \quad y^2 &\!\!=\!\!& (x+bc)(x+ac)(x+ab)  \\
 &\!\!=\!\!& x^3+(ab+ac+bc)x^2+(a^2bc+ab^2c+abc^2)x+a^2b^2c^2 \nonumber
\end{eqnarray}
in the Weierstrass form. There are three rational points on $E'$ of order $2$:
$$ T_1=[-bc,0], \quad T_2=[-ac,0], \quad T_3=[-ab,0], $$
and also other obvious rational points
$$ P=[0,abc], \quad Q=[1,rst]. $$
It is easy to verify that $Q=2R$, where
 $$ R= [rs+rt+st+1, (r+s)(r+t)(s+t)]. $$
In general, we may expect that the points $P$ and $R$ will be two independent
points of infinite order, and therefore that ${\rm rank}\,E(\mathbb{Q})\geq 2$.
Thus, assuming various standard conjectures, we may expect that the most of
elliptic curves induced by Diophantine triples with the above construction will have
the Mordell-Weil group $E(\mathbb{Q})$ isomorphic to
$\mathbb{Z}/2\mathbb{Z} \times \mathbb{Z}/2\mathbb{Z} \times {\mathbb{Z}}^2$ or
$\mathbb{Z}/2\mathbb{Z} \times \mathbb{Z}/2\mathbb{Z} \times {\mathbb{Z}}^3$.

The main purpose of this paper is to study which other groups are possible here.
Namely, we will investigate situations in which $P$ or $R$ have finite order,
or they are dependent, and in particular we would try to construct curves
with more independent points of infinite order.

According to Mazur's theorem, for the torsion group $E( \mathbb{Q})_{\rm tors}$ we have at most
four possibilities: $\mathbb{Z}/2\mathbb{Z} \times \mathbb{Z}/2m\mathbb{Z}$ with $m=1,2,3,4$.
In \cite{D-ja99}, it was shown that if $a,b,c$ are positive integers, then the cases
$m=2$ and $m=4$ are not possible. However, in the present paper we will show that
for $a,b,c$ non-zero rationals all four groups are indeed possible.

Let us also note that every Diophantine pair $\{a,b\}$ can be extended to
a Diophantine triple by a very simple extension: $c=a+b+2r$.
This construction was known already to Euler (and maybe even to Diophantus).
The direct computation shows that for triples of the form
$\{a,b,a+b+2r\}$, the points $P$ and $R$ are not independent, since $2P=-2R$.

\section{Search for elliptic curves with high rank}
\label{rank}

In last few years, several authors considered the problem of
construction of elliptic curves with some prescribed property and relatively high rank.
This includes curves with given torsion group (see \cite{D-table-high,Kul-Sta} and the references given there),
curves $y^2=x^3-n^2x$ related to congruent numbers \cite{Rog},
curves of the form $y^2=x^3+dx$ \cite{Elk-web},
curves $x^3+y^3=m$ related to the taxicab problem \cite{Elk-Rog},
curves $y^2=(ax+1)(bx+1)(cx+1)(dx+1)$ induced by Diophantine quadruples $\{a,b,c,d\}$ \cite{D-japan}, etc.

\bigskip

Let $G$ be an admissible torsion group for an elliptic curve over the rationals (according to Mazur's theorem).
Let us define
$$ B(G)=\sup \{ {\rm
rank}\,(E(\mathbb{Q})) \,:\, E( \mathbb{Q})_{\rm tors} = G
\}.
$$
The conjecture is that $B(G)$ is unbounded for all $G$. In the following table we give the best known lower bounds for $B(G)$.
Details on the record curves appearing in the table and full list of the references can be found on the
author's web page \cite{D-table-high}.

\bigskip

\begin{center}
{\footnotesize

\begin{tabular}{|@{\quad}c@{\quad}|@{\quad}c@{\quad}|@{\quad}l@{\quad}|}
\hline\rule{0pt}{12pt} $G$ & $B(G) \geq $ &
Author(s)\\[2pt]
\hline\rule{0pt}{12pt}
0 &          28 &        Elkies (2006) \\
  $\mathbb{Z}/2\mathbb{Z}$  &       18 &       Elkies  (2006) \\
  $\mathbb{Z}/3\mathbb{Z}$  &       13 &      Eroshkin  (2007) \\
  $\mathbb{Z}/4\mathbb{Z}$ &       12 &     Elkies (2006)  \\
  $\mathbb{Z}/5\mathbb{Z}$  &        6 &       Dujella \& Lecacheux  (2001) \\
  $\mathbb{Z}/6\mathbb{Z}$  &        7 &       Dujella  (2001,2006) \\
  $\mathbb{Z}/7\mathbb{Z}$  &        5 &       Dujella \& Kulesz (2001), Elkies (2006) \\
  $\mathbb{Z}/8\mathbb{Z}$  &        6 &       Elkies  (2006) \\
  $\mathbb{Z}/9\mathbb{Z}$  &        3 &       Dujella (2001), MacLeod (2004), Eroshkin (2006) \\
  $\mathbb{Z}/10\mathbb{Z}$ &        4 &       Dujella (2005), Elkies (2006) \\
  $\mathbb{Z}/12\mathbb{Z}$ &        3 &       Dujella (2001,2005,2006), Rathbun (2003,2006) \\
$\mathbb{Z}/2\mathbb{Z} \times \mathbb{Z}/2\mathbb{Z}$ &    14 &       Elkies  (2005) \\
$\mathbb{Z}/2\mathbb{Z} \times \mathbb{Z}/4\mathbb{Z}$ &     8 &       Elkies  (2005) \\
$\mathbb{Z}/2\mathbb{Z} \times \mathbb{Z}/6\mathbb{Z}$ &     6 &       Elkies (2006) \\
$\mathbb{Z}/2\mathbb{Z} \times \mathbb{Z}/8\mathbb{Z}$ &     3 &       Connell (2000), Dujella (2000,2001,2006), \\
  & & Campbell \& Goins (2003), Rathbun (2003,2006) \\[2pt]
\hline
\end{tabular}
}

\end{center}
\bigskip

We will now briefly describe the main steps in the construction of high rank curves
with prescribed properties. These steps have been already applied, with various
modifications, in obtaining curves from the above table, and we will also apply them
in the following sections.

The first step is to find a parametric family of elliptic curves over $\mathbb{Q}$
which contains curves with relatively high rank
(i.e. an elliptic curve over the field of rational functions
$\mathbb{Q}(T)$ with large generic rank) which satisfy the prescribed property. Here it is not
always the best idea to use the family with largest known generic rank, since these families
usually contain curves with very large coefficients for which it is very hard to compute the rank.

In the second step we want to find, in the given family of curves, the best candidates for higher rank.
The main idea here is that a curve is more likely to have large rank if $\# E(\mathbb{F}_p)$ is
relatively large for many primes $p$. We will use the following realization od this idea.
For a prime $p$, we put $a_p=a_p(E)=p+1- \# E(\mathbb{F}_p)$. For a fixed
integer $N$, we define
\begin{eqnarray*}
S(N,E) = \!\!\! \sum_{p\leq N,~p ~{\rm prime}} \left(1-\frac{p-1}{
\#E(\mathbb{F}_p) } \right) \log(p)= \!\!\! \sum_{p\le N,~p ~{\rm
prime}} \frac{-a_p+2}{p+1-a_p} ~\log(p).
\end{eqnarray*}
It is experimentally known (see \cite{Mestre,Nagao,D-r15}) that
we may expect that high rank curves have large $S(N,E)$. In
\cite{Camp}, some arguments were given which show that the
Birch and Swinnerton-Dyer conjecture gives support to this observation.
The sum $S(N,E)$ can be very efficiently computed (e.g. using
PARI \cite{PARI}) for $N<10000$. After this sieving method, we may
continue to investigate the best, let us say, $1 \%$ of curves. Since, we are working
with curves with torsion points of order $2$, we may compute the Selmer rank
for these curves, which is well-known upper bound for the actual rank of the curve.
This can be done using an appropriate option in Cremona's program {\it mwrank} \cite{Cremona}.

Only for the curves for which that upper bound is satisfactory large, we try to
compute the rank exactly. Again, the best available software for this purpose
is {\it mwrank} which uses 2-descent (via 2-isogeny if possible) to
determine the rank, obtain a set of
points which generate $E(\mathbb{Q})$ modulo $2E(\mathbb{Q})$, and finally saturate to a
full $\mathbb{Z}$-basis for $E(\mathbb{Q})$. The program package {\it APECS} \cite{APECS} and a program
that implements LLL reduction on the lattice of points of $E$,
provided by Rathbun \cite{Rath}, is used to reduce the heights of the generators.
In the cases when 4-descent is appropriate to perform
(for curves with a torsion point of order $4$, and with a generator
of very large height) we used an implementation of 4-descent in {\it MAGMA} \cite{MAGMA}.

\section{Torsion group $\mathbb{Z}/2\mathbb{Z} \times \mathbb{Z}/2\mathbb{Z}$}
\label{Z2*Z2}

In \cite{D-r7}, we have constructed a parametric family of elliptic curve with the
torsion group isomorphic to $\mathbb{Z}/2\mathbb{Z} \times \mathbb{Z}/2\mathbb{Z}$
and the generic rank $\geq 4$. The construction started with a Diophantine triple
$\{a,b,c\}$. We assigned to this triple the elliptic curve $E'$ as in (\ref{E'}), and
defined $d=x(P+Q)$, $e=x(P-Q)$. If $d,e\neq 0$, then $\{a,b,c,d\}$ and $\{a,b,c,e\}$ are Diophantine
quadruples (see \cite{D-ja99}). If $ed+1$ is a perfect square (and in \cite{D-r7} a
parametric solution to this equation was found), then we may expect that the elliptic curve
$y^2 = (bx + 1)(dx + 1)(ex + 1)$, induced by the Diophantine triple $\{b,d,e\}$,
has at least four independent points of infinite order, namely, points with $x$-coordinates
$0$, $a$, $c$ and $1/(bde)$. By a specialization, we found an elliptic curve of rank $7$
in that family. Here we will improve that result and construct several elliptic curves
of the form (\ref{ell}) which have rank equal to $8$ or $9$.

\medskip

The well-known family of Diophantine quadruples
\begin{equation} \label{kkk}
 \{ k-1, k+1, 4k, 16k^3-4k \}, \quad (k\in \mathbb{Z}, \,\, k\geq 2)
\end{equation}
has been studied by several authors. In \cite{D-deb}, it was proved that
the fourth element in this quadruple in unique, i.e. if $\{k-1,k+1,4k,d\}$
is an (integer) Diophantine quadruple, then $d=16k^3-4k$ (see also \cite{Fuj,BDM}).
It seems natural to consider the families of elliptic curves induced
by (\ref{kkk}). However, in \cite{D-param} it was shown that the triple
$\{ k-1, k+1, 4k \}$ induces an elliptic curve over $\mathbb{Q}(k)$ with generic rank
equal to $1$ (this agrees with the fact that $\{ k-1, k+1, 4k \}$ is of the form
$\{ a, b, a+b+2r \}$). Therefore, we will try to obtain curves with the higher rank induced by
other subtriples of (\ref{kkk}).

We first consider the family of Diophantine triples
$$  \{ k-1, k+1, 16k^3-4k \},  \quad (k\in \mathbb{Q}).$$
Applying the strategy described in Section \ref{rank}, we find a curve
with the rank equal to $9$ for $k=3593/2323$.
We have the Diophantine triple
$$ \left\{ \frac{1270}{2323} , \,\, \frac{5916}{2323}, \,\, \frac{664593861324}{12535672267} \right\}, $$
and the corresponding elliptic curve (its minimal Weierstrass equation)
\begin{eqnarray*}
y^2 &=& x^3 - 263759257625979218346701293692x
\\ & & \mbox{}+ 43309770676275925610968063087567021709640976.
\end{eqnarray*}
Torsion points are
\begin{eqnarray*}
& \mathcal{O}, [391223566189142, 0], [-581574668058484, 0], [190351101869342, 0],
\end{eqnarray*}
while independent points of infinite order are
\begin{eqnarray*}
& [13382356740992, 6307332780932304905700], \\
& [137392393772492, 3108820636640783206800], \\
&  [151121694899342, 2627033807399227434000], \\
&  [182390949979797, 1126929502996358494505], \\
&  [-100285963891570, 8291713851182161095696], \\
&  [638038681834022, 11608723965551290530480], \\
&  [-570129376204450, 2892670061337006977376], \\
&  [-581493416436883, 246987048216416159925], \\
&  [395944953729830, 974097947650374250704].
\end{eqnarray*}

We also found several examples with the rank equal to $8$ in this family:
$k=286/69$, $69/1144$, $1169/1268$, $1225/1959$, $1443/1156$, $1981/1941$,
$2447/50$, $4350/1159$, $5781/782$.

Next we consider the family
$$ \{ k-1, 4k, 16k^3-4k \}, \quad (k\in \mathbb{Q}). $$
(Note that the triples of the form $ \{ k+1, 4k, 16k^3-4k \}$ induce
the same family, by the correspondence $k \leftrightarrow -k$.)
In this case, we find a curve
with the rank equal to $9$ for $k=-2673/491$,
and several examples with the rank equal to $8$, e.g. for
$k=65/521$,\,
$864/1415$,\, $909/2741$,\, $1500/2339$,\, $1610/4401$,\,
$1914/2969$,\,
$3656/5127$,\,
$4435/3378$,\, 
$6648/3473$,\,
$-175/2098$,\,
$-291/674$,\,
$-338/911$, \linebreak
$-470/889$,\, $-535/5178$,\, $-559/807$,\, 
$-705/1703$,\,  
$-1224/4555$,\,
$-1443/964$,\,
$-1610/1629$,\, 
$-2123/4703$,\,
$-2209/2927$.

\medskip

We may also consider rational Diophantine triples of the form $\{1,3,c\}$.
Here we find two examples with rank equal to $8$ for
$$ c=\frac{5043716589720}{9928996362961} \quad \mbox{and}\quad
c=\frac{507857302680}{1680262081}.$$

\medskip

Let us mention that Gibbs discovered 46 examples of rational
Diophantine sextuples (\cite{Gib1,Gib2}). By computing the ranks
for all curves of the form $y^2=(ax+1)(bx+1)(cx+1)$, where
$\{a,b,c\}$ is a subtriple of some of the Gibbs's sextuples, we
find a curve of rank $8$ for the Diophantine triple
$$ \left\{ \frac{494}{35},\,\, \frac{1254396}{665},\,\, \frac{11451300}{5067001} \right\}. $$

Of course, the curves with rank less than $7$
(and greater than $0$)
are easy to find, and they already appeared in the literature (see \cite{D-param,D-P,D-r7}).
Therefore, we can summarize the results from this section in the following proposition:

\begin{proposition} \label{z2z2}
For each $1\leq r \leq 9$, there exists a Diophantine triple $\{a,b,c\}$
such that the elliptic curve $y^2=(ax+1)(bx+1)(cx+1)$ has the torsion group isomorphic to
$\mathbb{Z}/2\mathbb{Z} \times \mathbb{Z}/2\mathbb{Z}$ and the rank equal to $r$.
\end{proposition}

\section{Torsion group $\mathbb{Z}/2\mathbb{Z} \times \mathbb{Z}/4\mathbb{Z}$}
\label{Z2*Z4}

In this section, we consider elliptic curves with a torsion subgroup isomorphic
to $\mathbb{Z}/2\mathbb{Z} \times \mathbb{Z}/4\mathbb{Z}$. It follows from the 2-descent proposition
(see \cite[4.1, p.37]{Hus}, \cite[4.2, p.85]{Kna}), that such curve has the
equation of the form
\begin{equation} \label{eq:z2z4}
y^2=x(x+x_1^2)(x+x_2^2), \quad x_1,x_2 \in \mathbb{Q}.
\end{equation}
The point $[0,0]$ is a double point (i.e. point of the form $2S$, where $S$ is a rational point on the curve) of order $2$.
Translating the elliptic curve (\ref{E'}) induced by the Diophantine triple $\{a,b,c\}$,
we obtain the equation
\begin{equation} \label{eq:abc24}
y^2=x(x+ac-ab)(x+bc-ab).
\end{equation}
Therefore, if we can find $a,b,c$ such that $ac-ab$ and $bc-ab$ are perfect squares, then the elliptic curve
induced by $\{a,b,c\}$ will have a torsion subgroup isomorphic
to $\mathbb{Z}/2\mathbb{Z} \times \mathbb{Z}/4\mathbb{Z}$.

A simple way to fulfill these conditions is to choose $a$ and $b$ such that $ab=-1$.
Then $ac-ab=ac+1=s^2$ and $bc-ab=bc+1=t^2$. It remains to find $c$ such that
$\{a, -1/a, c\}$ is a Diophantine triple. Using the standard extension $c=a+b+2r$,
we may take $c=a-\frac{1}{a}$. However, it is easy to prove, using Shioda's formula (\cite{Shi}),
that the family of elliptic curves
$$ y^2=x^3+ (a^4+1)x^2+a^4x, $$
obtained with this construction, has the generic rank equal to $0$.
We may ask what happened with the points $P=[0,abc]$ and $Q=[1,rst]$.
It is easy to see that $2P=Q=T_3$. Hence, in this case $P$ and $Q$ are indeed points of finite order.
We are able to find
examples with rank equal to $0,1,2,3,4$ is this family, but in order to find curves with
higher rank, we will consider some other constructions.

We are searching for parametric solutions of the system
\begin{equation} \label{eq:ac}
 ac+1 = \square, \quad -\frac{c}{a} +1=\square.
\end{equation}
Multiplying these two conditions, we obtain
\begin{equation} \label{eq:ace}
 a(ac+1)(a-c)= \square,
\end{equation}
which, for given $c$, may be regarded as an elliptic curve. We already know one (non-torsion)
parametric solution of (\ref{eq:ace}), namely $a=T$, $c=T-\frac{1}{T}$. By duplicating the
corresponding point on the elliptic curve (\ref{eq:ace}), we obtain another solution
$a=\frac{(T^2+1)^2(T^2-1)}{4T^3}$, with the same $c$. By the 2-descent proposition,
these values of $a$ and $c$ also satisfy the original system (\ref{eq:ac}).
We have again that $Q=T_3$, but now the point $P$ has infinite order.

By searching for curves with high rank in this family of elliptic curves, with the methods described
in Section \ref{rank}, we are able to
find two curves with rank equal to 5, for $T=12/5$ and $T=24/7$, corresponding to
the Diophantine triples
$$ \left\{ \frac{3398759}{864000}, \,\, -\frac{864000}{2298759}, \,\, \frac{119}{60} \right\}, $$
$$ \left\{ \frac{205859375}{18966528}, \,\, -\frac{18966528}{205859375}, \,\, \frac{527}{168} \right\}.$$

We will improve this result by considering a different parametric solution of the system (\ref{eq:ac}).
Inserting $ac+1=s^2$ into $-\frac{c}{a}+1=t^2$, we obtain
$$ 1-s^2+a^2=\square. $$
which has the solution of the form
$$ a=\frac{\alpha T +1}{T - \alpha}, \quad s=\frac{T+ \alpha}{T - \alpha}. $$
We take $\alpha=2$, which gives
$$ a=\frac{2T+1}{T-2}, \quad b=\frac{2-T}{2T+1}, \quad
c=\frac{8T}{(2T+1)(T-2)}. $$
This yields again the family of elliptic curves with generic rank $\geq 1$. We are able to find in this family
several examples of curves with rank equal to 7, e.g. for $T=7995/6562$, $28853/5306$, $55204/28537$, $87046/1523$,
$95827/81626$, $134726/16613$, and many examples with rank equal to $6$,
e.g. for $T=399/160$, $452/173$, $698/561$, $1212/661$, $1253/974$,
$1263/707$, $1463/1081$.

We give some details only for $T=7995/6562$. In that case we obtain the Diophantine triple
$$ \left\{ \frac{22552}{5129} , \,\, -\frac{5129}{22552}, \,\, \frac{52463190}{14458651} \right\}, $$
and the corresponding elliptic curve (its minimal Weierstrass equation) is
\begin{eqnarray*}
y^2 &=& x^3 - x^2 - 66316100370037243788808101431860x
\\ & & \mbox{}+ 207787397329581777063539110158423853553882263492.
\end{eqnarray*}
Torsion points:
\begin{eqnarray*}
& \mathcal{O}, [4628372593789489, 0], [4774517796226298, 0], [-9402890390015786, 0], \\
& [6213948304937164, 188792790675737045056350], \\
& [6213948304937164, -188792790675737045056350],  \\
& [3335087287515432, -153990087256491803556554], \\
& [3335087287515432, 153990087256491803556554].
\end{eqnarray*}
Independent points of infinite order:
\begin{eqnarray*}
& [2296719982411009, 259986354975552263423820], \\
& [-754552913305211, 507342881379999424501350], \\
& [4623547820142394, 3196373093851098851280], \\
& [4774647845487904, 519325665120008067090], \\
& [2882111728866944, 201488073470229552730830], \\
& [4572248586935254, 12595554476372830774560], \\
& [-14146681940011859/4, 5047542184921547691312855/8].
\end{eqnarray*}

Hence, we have proved the following result.

\begin{proposition} \label{z2z4}
For each $0\leq r \leq 7$, there exists a Diophantine triple $\{a,b,c\}$
such that the elliptic curve $y^2=(ax+1)(bx+1)(cx+1)$ has the torsion group isomorphic to
$\mathbb{Z}/2\mathbb{Z} \times \mathbb{Z}/4\mathbb{Z}$ and the rank equal to $r$.
\end{proposition}

\section{Torsion group $\mathbb{Z}/2\mathbb{Z} \times \mathbb{Z}/6\mathbb{Z}$}
\label{Z2*Z6}

General form of curves with the torsion group isomorphic to $\mathbb{Z}/2\mathbb{Z} \times \mathbb{Z}/6\mathbb{Z}$
is
$$ y^2= (x+\alpha^2)(x+\beta^2)\left(x+\frac{\alpha^2\beta^2}{(\alpha-\beta)^2}\right) $$
(then the point $[0,\alpha^2\beta^2/(\alpha-\beta)]$ is of order 3; see \cite{Kul-acta}).
Let us force a Diophantine triple equation
\begin{equation} \label{XY26}
 y^2=(x+bc)(x+ac)(x+ab)
\end{equation}
to have this form.

Comparison gives
\begin{align}
& \alpha^2+1=bc+1=t^2, \label{alpha} \\
& \beta^2+1=ac+1=s^2, \label{beta} \\
& \alpha^2\beta^2+(\alpha-\beta)^2 = \square. \label{treca}
\end{align}
The first two conditions (\ref{alpha}) and (\ref{beta}) have parametric
solutions
$$ \alpha=\frac{2u}{u^2-1}, \quad \beta=\frac{v^2-1}{2v}.$$
Inserting this
into (\ref{treca}), we obtain the equation $F(u,v)=z^2$, where
\begin{eqnarray} \label{Fuv}
  \,\,\,\,F(u,v) &\!\!\!=\!\!\!& (v^4-2v^2+1)u^4+(-8v^3+8v)u^3 \\
& &\mbox{}+ (2v^4+2+12v^2)u^2+ (-8v+8v^3)u+v^4-2v^2+1. \nonumber
\end{eqnarray}

The condition $F(u,v)=z^2$ is satisfied e.g. for
$$ u=\frac{v^3+v}{v^2-1}. $$
Indeed,
then $F(u,v)=\frac{(v^6-v^4+3v^2+1)^2}{(v-1)^2 (v+1)^2}$.
Hence, for
$$ \alpha=\frac{2T^5-2T}{T^6+T^4+3T^2-1}, \quad
\beta=\frac{T^2-1}{2T}, $$
we obtain the parametric family of Diophantine triple equations with
the torsion group isomorphic to $\mathbb{Z}/2\mathbb{Z} \times \mathbb{Z}/6\mathbb{Z}$.
It is easy to check that the point $[1,rst]$ on (\ref{XY26}) has an infinite order
(by finding a suitable specialization, or by listing explicitly all 12 torsion points
on (\ref{XY26}) over $\mathbb{Q}(T)$).
Hence, we found an elliptic curve over $\mathbb{Q}(T)$ with the torsion group
$\mathbb{Z}/2\mathbb{Z} \times \mathbb{Z}/6\mathbb{Z}$ and the generic rank $\geq 1$,
which ties the current record (\cite{Camp,Kul-acta,Lec2628}).

Since the constructed curve over $\mathbb{Q}(T)$ has very large coefficients, it is
not surprising that we are able to compute the rank only for few specializations, and
among them we find examples with rank equal to $1,2,3$. Rank $3$ is obtained for $T=7$,
which corresponds for the Diophantine triple
$$ \left\{ \frac{721176}{193193}, \,\, \frac{20580000}{829322351},
\,\, \frac{662376}{210343} \right\}. $$

Instead of using a parametric solution, we can also try to search for solutions $u,v$
of the equation (\ref{Fuv}) with small numerators and denominators, and to compute the rank of corresponding elliptic curves.
Using this approach we are able to find a curve with rank equal to $4$. It is obtained for
$u=34/35$ and $v=8$, i.e. for the Diophantine triple
$$ \left\{ \frac{39123}{96976}, \,\, \frac{12947200}{418209}, \,\,  \frac{42427}{1104} \right\}. $$
The curve is
\begin{eqnarray*}
y^2 + xy &=& x^3  - 24046649084795243589952562390x
\\ & & \mbox{}+ 1435226116741326558309046453105518735800100.
\end{eqnarray*}
Torsion points:
\begin{eqnarray*}
& \mathcal{O}, [-179058763357620, 89529381678810], \\
& [89873668514380, -44936834257190], \\
& [356740379372959/4, -356740379372959/8], \\
& [92726794888780, -52405873597247415590], \\
& [92726794888780, 52405780870452526810], \\
& [86369148214060, 51179899438633016410], \\
& [154777835944300, 1192150615832496114010], \\
& [37033400507980, -771678209256671722790], \\
& [86369148214060, -51179985807781230470], \\
& [154777835944300, -1192150770610332058310], \\
& [37033400507980, 771678172223271214810].
\end{eqnarray*}
Independent points of infinite order:
\begin{eqnarray*}
& [35387651068492, 792834860571692154586], \\
& [-39964997451020, -1527225651415581670190], \\
& [-8547561811220, -1280680222922667973190], \\
& [90070190194252, 6841914086525854426].
\end{eqnarray*}

\begin{proposition} \label{z2z6}
For each $1\leq r \leq 4$, there exists a Diophantine triple $\{a,b,c\}$
such that the elliptic curve $y^2=(ax+1)(bx+1)(cx+1)$ has the torsion group isomorphic to
$\mathbb{Z}/2\mathbb{Z} \times \mathbb{Z}/6\mathbb{Z}$ and the rank equal to $r$.
\end{proposition}

\section{Torsion group $\mathbb{Z}/2\mathbb{Z} \times \mathbb{Z}/8\mathbb{Z}$}
\label{Z2*Z8}

Finally, we consider the largest possible torsion group
$\mathbb{Z}/2\mathbb{Z} \times \mathbb{Z}/8\mathbb{Z}$. As we have already noted
in Section \ref{Z2*Z8}, the torsion group of elliptic curves induced by
Diophantine triples of the form
$$ \{a, -\frac{1}{a}, a-\frac{1}{a} \} $$
contains a subgroup isomorphic to $\mathbb{Z}/2\mathbb{Z} \times \mathbb{Z}/4\mathbb{Z}$.
In that case, the points of order $4$ on
$$ y^2=(x+ab)(x+ac)(x+bc) $$
are $P=[0,abc]$, $P+T_1$, $P+T_2$, $P+T_3$, where $T_1=[-bc,0]$, $T_2=[-ac,0]$,
$T_3=[-ab,0]$. Hence, our elliptic curve will have the torsion group isomorphic to
$\mathbb{Z}/2\mathbb{Z} \times \mathbb{Z}/8\mathbb{Z}$ if some of the points
$P,P+T_1,P+T_2,P+T_3$ is a double point. We will use 2-descent proposition again.
Consider the point $P+T_2$. It will be a double point iff $(b-a)(b-c)$ and $b(b-a)$ are
both perfect squares. These conditions lead to a single condition that $a^2+1$ is
a perfect square. Therefore, we have proved that all Diophantine triples of the form
$$ \left\{ \frac{2T}{T^2-1}, \,\, -\frac{1-T^2}{2T}, \,\, \frac{6T^2-T^4-1}{2T(T^2-1)} \right\}, \quad t\in \mathbb{Q}, $$
induce elliptic curves with torsion group isomorphic to
$\mathbb{Z}/2\mathbb{Z} \times \mathbb{Z}/8\mathbb{Z}$. The induced curves have the equation of the form
\begin{equation} \label{z2z8kul}
 y^2=x(x+s^2)(x+t^2)= x \left(x+ \left(\frac{2T}{T^2-1} \right)^2 \right)
 \left(x+ \left(\frac{T^2-1}{2T} \right)^2 \right).
\end{equation}
But, according to \cite{Kul-acta}, \emph{every} elliptic curve over $\mathbb{Q}$ with
torsion group $\mathbb{Z}/2\mathbb{Z} \times \mathbb{Z}/8\mathbb{Z}$ has an equation of the
form (\ref{z2z8kul}). Therefore, every such curve is induced by a Diophantine triple.
This fact has been independently proved by Campbell and Goins in \cite{C-G}.

Thus, we are left with the question which ranks are possible for the elliptic curves with
torsion group $\mathbb{Z}/2\mathbb{Z} \times \mathbb{Z}/8\mathbb{Z}$. It is known that there
exist infinitely many such curves with rank $\geq 1$ (see \cite{A-M,C-G,Kul-acta,Lec2628}), although no
such parametric family (curve over $\mathbb{Q}(T)$) is known. It is easy to find examples with
rank equal to $0,1,2$. The first example with rank equal to $3$ was found in 2000, independently,
by Connell \cite{APECS} and the author \cite{D-table-high}. It was the curve
\begin{eqnarray*}
y^2 + xy &=& x^3  - 15745932530829089880x
\\ & & \mbox{}+ 24028219957095969426339278400,
\end{eqnarray*}
with torsion points:
\begin{eqnarray*}
& \mathcal{O}, [-4581539664, 2290769832], [-1236230160, 203972501847720], \\
& [2132310660, 12167787556920], [2452514160, 12747996298920], \\
& [9535415580, 860741285907000], [2132310660, -12169919867580], \\
& [-1236230160, -203971265617560], [9535415580, -860750821322580], \\
& [2452514160, -12750448813080], [2346026160, -1173013080], \\
& [1471049760, 63627110794920], [1471049760, -63628581844680], \\
& [3221002560, -82025835631080], [3221002560, 82022614628520], \\
& [8942054015/4, -8942054015/8],
\end{eqnarray*}
and independent points of infinite order:
\begin{eqnarray*}
& [2188064030, -7124272297330], \\
& [396546810000/169, 1222553114825160/2197], \\
& [16652415739760/3481, 49537578975823615480/205379].
\end{eqnarray*}
The curve was induced by the Diophantine triple
$$ \left\{ \frac{408}{145}, \,\, -\frac{145}{408}, \,\, -\frac{145439}{59160} \right\}. $$

In the meantime, several other curves with rank equal to 3 were found by Rathbun \cite{Rath},
Campbell and Goins \cite{C-G}, and the author (see also \cite{BFHHZ}). Here we will mention our findings.
Using a similar search procedure, as in the previous sections, we have discovered elliptic curves
with the torsion group $\mathbb{Z}/2\mathbb{Z} \times \mathbb{Z}/8\mathbb{Z}$ and the rank
equal to $3$, which correspond to the following Diophantine triples:
$$ \left\{ \frac{1692}{1885}, \,\, -\frac{1885}{1692}, \,\, -\frac{690361}{3189420} \right\}, $$
$$ \left\{ \frac{79040}{35409}, \,\, -\frac{35409}{79040}, \,\, \frac{4993524319}{2798727360} \right\}, $$
$$ \left\{ \frac{77556}{59917}, \,\, -\frac{59917}{77556}, \,\, \frac{2424886247}{4646922852} \right\}, $$
$$ \left\{ \frac{128760}{176111}, \,\, -\frac{176111}{128760}, \,\, -\frac{14435946721}{22676052360} \right\}, $$
$$ \left\{ \frac{424580}{799029}, \,\, -\frac{799029}{424580}, \,\, -\frac{458179166441}{339251732820} \right\}, $$
$$ \left\{ \frac{451352}{974415}, \,\, -\frac{974415}{451352}, \,\, -\frac{745765964321}{439804159080} \right\}. $$
We give some details on the curve corresponding to the last triple, since this case
is the most technically involved and time consuming.
Its minimal Weierstrass equation is
{\small
\begin{eqnarray*}
& y^2 + xy = x^3  - 16188503722614063108729139735755154904562292360x \\
  \mbox{} & \!\!\! + \,786863421808206463969913495490892469346874709447053592901366525761600.
\end{eqnarray*} }
Torsion points are:
\begin{eqnarray*}
& \mathcal{O},  [126942113771663398101920, 27882680574001240245704236363397240], \\
& [126942113771663398101920, -27882680574128182359475899761499160], \\
& [30228264599630424878720, -18031474247759343557251945104801160], \\
& [30228264599630424878720, 18031474247729115292652314679922440], \\
& [422083239655931288586320, -262963676295121325354530570456161160], \\
& [85392774125986994678180, -5211521754389769127545451481401720], \\
& [-61613633858042970798640, 39375066083178460970956312648094840], \\
& [61906075162778279941820, 4684356135071180542946299694504840], \\
& [85392774125986994678180, 5211521754304376353419464486723540], \\
& [422083239655931288586320, 262963676294699242114874639167574840], \\
& [61906075162778279941820, -4684356135133086618109077974446660], \\
& [-61613633858042970798640, -39375066083116847337098269677296200], \\
& [68209869346874809632016, -34104934673437404816008], \\
& [78585189185646911490320, -39292594592823455745160], \\
&[-587180234130086884489345/4, 587180234130086884489345/8],
\end{eqnarray*}
while independent points of infinite order are:
{\small
\begin{eqnarray*}
P_1 &\!\!=\!\!& [66119657073815781066800, -2355339128918565969721076384104840], \\
P_2 &\!\!=\!\!& [401169287265672834550867500080/76335169, \\
&&      -17669918374394464360754418260810255353172665480/666940371553],
\end{eqnarray*} } {\tiny
\begin{eqnarray*}
P_3 &\!\!\!\!=\!\!\!\!& [2896876482219215018082911215320035879728511402808... \\
&& ...52543070142236950355376437227490536812623183120/ \\
&&      17156914675194799164872812895696296908615977178623591714662702334872955881, \\
&&      -1617954908834725344870195603458416032731793227363235820394569809657744633819... \\
&& ...656062208905822304087388210043858968407239618913517436708870237508440/ \\
&&      7106549528400313535484236316357837399782369803809731948... \\
&& ...0954479040616074792502966798554498341539450071275232779].
\end{eqnarray*} }%
(Numerators and denominators too large to fit in one line are given in two lines, which is indicated by $...$ at the end of
the first line and beginning of the second line.)
In this case, we were not able to compute the exact rank using {\it mwrank}
(we obtained that rank is equal to $2$ or $3$).
Namely, the coordinates of the points $P_3$ are too large to be found by 2-descent.
Therefore, here we used 4-descent implemented in {\it MAGMA}.

\medskip

Let us summarize the results from this section.

\begin{proposition} \label{z2z8}
For each $0\leq r \leq 3$, there exists a Diophantine triple $\{a,b,c\}$
such that the elliptic curve $y^2=(ax+1)(bx+1)(cx+1)$ has the torsion group isomorphic to
$\mathbb{Z}/2\mathbb{Z} \times \mathbb{Z}/8\mathbb{Z}$ and the rank equal to $r$.
\end{proposition}

{\bf Acknowledgements.} The author was supported by the Ministry
of Science, Education and Sports, Republic of Croatia, grants
0037110 and 037-0372781-2821. The experiments reported in this
paper were performed on the computers of the Laboratory for
Advanced Computations, Department of Mathematics, University of
Zagreb. The author would like to thank to Ivica Gusi\'{c} for many
useful comments on the first version of the manuscript and to
Randall Rathbun for his generous help with computations involving
4-descent, which were essential in computing the exact rank for
two curves of rank equal to $3$ in Section \ref{Z2*Z8}.

\bigskip

{\small \noindent
Andrej Dujella \\
Department of Mathematics \\ University of Zagreb
\\ Bijeni\v{c}ka cesta 30 \\
10000 Zagreb, Croatia \\
{\em E-mail address}: {\tt duje@math.hr}}

\end{document}